\documentclass[12pt]{article}
   \def\sqr#1#2{$\vcenter{\hrule height.#2pt
   \hbox{\vrule width.#2pt height#1pt \kern#1pt
   \vrule width.#2pt}
   \hrule height.#2pt}$}
   
   \newfont{\eightrm}{cmr8}
   \newfont{\ninerm}{cmr9}
   \newfont{\nineit}{cmti9}
   \newfont{\eightit}{cmti8}
 \def\no{\noindent}

\begin{document}
\begin{center}
\vskip 0.6cm {\large\bf Finding Short Cycles in an Embedded Graph in
Polynomial Time \footnote{Supported by NNSF of China under granted
number 10671073 and partially supported by Science and Technology
Commission of Shanghai Municipality (07XD14011)}}

\vskip 0.5cm

{Han Ren \footnote{Supported by Shanghai Leading Academic Discipline
Project, Project Number£ºB407}\quad and \quad Ni Cao}

\vskip 0.3cm

 {\footnotesize{\it Department of Mathematics,East
China Normal University,}}

\vskip 0.3cm

{\footnotesize{\it{Shanghai 200062,P.R.China}}}

\end{center}

\vspace{0.5in}
\begin{center}
{\eightrm ABSTRACT}
\end{center}

\medskip
\begin{center}
\begin{minipage}{120mm}

Let ${\cal{C}}_1$ be the set of fundamental cycles of
breadth-first-search trees in a graph $G$ and ${\cal{C}}_2$ the set
of the sums of two cycles in ${\cal{C}}_1$. Then we show that
$(1)\,{\cal{C}}={\cal{C}}_1\bigcup{\cal{C}}_2$ contains a shortest
$\Pi$-twosided cycle in a $\Pi$-embedded graph $G$;$(2)$\, $\cal{C}$
contains all the possible shortest even cycles in a graph
$G$;$(3)$\,If a shortest cycle in a graph $G$ is an odd cycle, then
$\cal{C}$ contains all the shortest odd cycles in $G$. This implies
the existence of a polynomially bounded algorithm to find a shortest
$\Pi-$twosided cycle in an embedded graph and thus solves an open
problem of B.Mohar and C.Thomassen[2,pp112]

 \vskip 0.2cm

 \noindent{\rm\bf Key Words}\quad{$\Pi-$twosided cycle,
breadth-first-search tree, embedded graph.}

\vskip 0.2cm AMS {\rm\bf Classification}: (2000)05C10,05C30,05C45.
\end{minipage}
\end{center}


\newpage


\noindent{\bf 1.\quad Introduction}

\vskip 0.5cm

C.Thomassen showed that if cycles in a set of cycles satisfy the
{\it $3$-path-condition}, then there exists a polynomial time
algorithm that finds a shortest cycle in this set[3]. As
applications, he showed that the following types of shortest cycles
may be found in polynomial time.

\vskip 0.3cm {\it
\noindent (1)\quad A shortest $\Pi-$noncontractible cycle in a $\Pi-$embedded graph;\\
\noindent (2)\quad A shortest $\Pi-$nonseparating cycle in a $\Pi-$embedded graph;\\
\noindent (3)\quad A shortest $\Pi-$onesided cycle in a $\Pi-$
embedded graph. } \vskip 0.3cm

But what about a family of cycles which do not satisfy the
$3$-path-condition? In their monograph[2, pp112], B.Mohar and
C.Thomassen raised the following open problems:

\vskip 0.2cm \vskip 0.2cm {\it \noindent(a)\quad Is there a
polynomially bounded algorithm that finds
a shortest $\Pi-$contractible cycle in a $\Pi-$embedded graph?\\
\noindent(b)\quad Is there a polynomially bounded algorithm that
finds a shortest $\Pi-$surface-separating cycle in a $\Pi-$embedded graph?\\
\noindent(c)\quad Is there a polynomially bounded algorithm that
finds a shortest $\Pi-$twosided cycle in a $\Pi-$embedded graph? }

\vskip 0.3cm

Here in this paper we consider connected graphs and all the concepts
used are standard following from [1,2].

\vskip 0.2cm

Let $G$ be a $\Pi-$embedded graph and $T_x$ a {\it
breadth-first-search tree} rooted at a vertex $x$ of $G$. Let
$$
{\cal{C}}_1=\{C(T_x) | \forall x\in{V(G)}, C(T_x)\,\,is\,\, a\,\,
fundamental\,\,cycle\,\, of\,\, T_x\};
$$
$$
{\cal{C}}_2=\{ C | \exists x,y\in{V(G)}, C=C(T_x)\oplus C(T_y)\};\\
{\cal C}={\cal C}_1\bigcup{\cal C}_2,
$$
where "$\oplus$" is the operation defined as "$A\oplus
B=(A-B)\cup(B-A)$" for  any subsets $A,B$ of $E(G)$.

\vskip 0.2cm

\noindent{\bf Theorem A.}\quad {\it The collection $\cal C$ of
cycles defined above satisfies the following conditions:\\
\noindent (a). There exists a shortest $\Pi-$twosided cycle (in a
$\Pi-$embedded graph ) in $\cal C$; \\
\noindent (b). Each shortest even cycle of a graph is
contained in $\cal C$;\\
\noindent (c). If a shortest cycle of a graph is an odd cycle, then
every shortest odd cycle is contained in ${\cal C}_1$.}

\vskip 0.2cm

Therefore, Theorem A implies the existence of a polynomially bounded
algorithm to find short cycles defined above.

\vskip 0.2cm

\noindent{\bf Theorem B.}\quad{\it There exists a polynomially
bounded algorithm to find a shortest $\Pi-$twosided cycle in a
$\Pi-$embedded graph and all the shortest even cycles in a graph.}

 \vskip 0.2cm

This solves problem (c). Since every possible two-sided cycle in
an embedded graph in the projective plane is contractible, we have
the following

\vskip 0.2cm

\noindent{\bf Corollary.}\quad{\it There is a polynomially bounded
algorithm to find a shortest contractible cycle in an embedded graph
in the projective plane. }

\vskip 0.2cm

This answers the problems (a)-(c) in the case of projective plane
graphs.

\vskip 0.5cm

\noindent{\bf 2.\quad Proof of Main Result }

\vskip 0.3cm

A {\it generalized embedding scheme} of a graph $G$ in a surface
$\Sigma$ is a set of a rotation systems $\pi=\{{\pi}_v|v\in
{V(G)}\}$ together with a mapping $\lambda:
E(G)\rightarrow\{-1,+1\}$,called a {\it signature}, where ${\pi}_v$
is a clockwise ordering of edges incident with $v$. We define
$\Pi=(\pi,\lambda)$ as an embedding scheme for the graph $G$. A
cycle $C$ of an embedded graph $G$ ia called $\Pi-$twosided if $C$
contains even number of edges with negative signature;otherwise, it
is called $\Pi-$onesided.

\vskip 0.2cm

For a vertex $v\in{V(G)}$, we may change the clockwise ordering to
anticlockwise, i.e., ${\pi}_v$ is replaced by its inverse
${\pi}^{-1}_v$, and\,$\lambda(e)$ is replaced by $-\lambda(e)$ for
each edge $e$ incident with $v$. Therefore, we obtain another
embedding scheme ${\Pi}'=({\pi}',{\lambda}')$. It is clear that a
cycle is $\Pi-$twosided if and only if it is ${\Pi}'$-twosided. Two
of such embedding schemes are {\it equivalent} if and only if one
can be changed into another by a sequence of such {\it local
changes}. Thus, for any spanning tree $T$ of an embedded graph $G$,
we may always assume that each edge $e$ of $T$ has signature
$\lambda(e)=+1$ for convenience.

\vskip 0.2cm

The key of the proof of Theorem A is the following classification of
short cycles according to the distance between two vertices on such
cycles. Let $C$ be a shortest $\Pi-$twosided cycle in an embedded
graph $G$. Then $C$ must satisfy one of the following
conditions:\\
\noindent(1).$\forall x,y\in{V(C)}\Rightarrow d_C(x,y)=d_G(x,y);$\\
\noindent(2).$\exists x,y\in{V(C)}\Rightarrow d_C(x,y)>d_G(x,y$.

\vskip 0.2cm

\noindent{\bf Remark:}\,In the following, we will see that a cycle
satisfying (2) may be written as a sum of two shorter cycles.
Therefore, a shortest cycle in a collection of cycles satisfying the
$3$-path-condition can't satisfy (2).

\vskip 0.2cm

\noindent{\bf Lemma 1.}\quad{\it Let $G$ and $\cal C$ be as defined
in Theorem A and $C$ a shortest $\Pi$-twosided cycle in $G$
satisfying (2). Then $\cal C$ contains a shortest $\Pi-$twosided
cycle of $G$.}

\vskip 0.2cm

\noindent{\bf Proof of Lemma 1.}\quad We assume that $C$ has a
clockwise ( anticlockwise )orientation $\overrightarrow{C}$
(\overleftarrow{C})and for any two vertices $u$ and $v$ of $C$,
$u\overrightarrow{C}v$ ($u\overleftarrow{C}v$) denotes the closed
interval from $u$ to $v$ along $\overrightarrow{C}$
($\overleftarrow{C}$).  Similarly, we may define a segment $uPv$ as
the closed interval from $u$ to $v$ in a path $P$. If $u$ is a
vertex of $C$, then $u^{-} (u^{+})$ is the predecessor ( successor)
of $u$ along $\overrightarrow{C}$.

\vskip 0.2cm

We consider two vertices $x, y$ of $C$ such that $d_G(x,y)$ is
minimum subject to (2). Then for any shortest $(x-y)$ path $P$ in
$G$, $P$ and $C$ has no inner vertex on $C$( since otherwise, there
would be another two vertices $x_1, y_1$ od $C$ satisfying (2) and
$D_G(x_1,y_1)<d_G(x,y)$). Let $T_x$ be a breadth-first-search tree
rooted at $x$. Then it contains a $(x-y)$ path P and each of the two
cycles $P\cup x\overrightarrow{C}y$ and  $P\cup x\overleftarrow{C}y$
is a $\Pi-$onesided cycle shorter than that of $C$ and satisfies (1)
(since otherwise one of them may be written as a sum of two shorter
cycles,among then, one is a $\Pi-$twosided ). Under these
structure,each of $P\cup x\overrightarrow{C}y$ and  $P\cup
x\overleftarrow{C}y$ is a fundamental cycle of a
breadth-first-search tree rooted at a vertex on itself. This
completes the proof of Lemma 1.

 \vskip 0.2cm

 Now we turn to the cycles satisfying (1). Firstly, the following
 result is easy to be verified and we omit the proof of it.

 \vskip 0.2cm

 \noindent{\bf Lemma 2.}\quad{\it Let $C$ be a shortest $\Pi-$twosided cycle in a
 $\Pi-$embedded graph $G$. Then for any three vertices $x,u,v$ of $C$ with
$$
d_c(x,u)=d_c(x,v)=[\frac{|C|-1}{2}],
$$
any shortest $(x-u)$path $P$ and $x\overleftarrow{C}v$ has no inner
vertex in common.}

\vskip 0.2cm

Thus, for a shortest $\Pi-$twosided cycle $C$ satisfying the
conditions in Lemma 2 and the last common vertex $\alpha$ ( called
{\it branched vertex}) contained in $(x-u)$path $P_1$ and
$(x-v)$path $P_2$ of $T_x$,
$x\overrightarrow{P_1}\alpha(=x\overrightarrow{P_2}\alpha)$ has no
inner vertex on $x\overleftarrow{C}v( x\overrightarrow{C}u)$.

\vskip 0.2cm

\noindent{\bf Lemma 3.}\quad{\it Let $G$ and $\cal C$ be as defined
in Theorem A and $C$ a shortest $\Pi-$twosided cycle in $G$
satisfying (1). Then $\cal C$ contains a shortest $\Pi-twosided$
cycle of $G$.}

\vskip 0.2cm

\newpage

\noindent{\bf Proof of Lemma 3.}

\vskip 0.2cm

\noindent{\bf Case 1}\quad $|C|=2n,\, n\in N$.

\vskip 0.2cm

 Let $x$ and $y$ be
two vertices of $C$ with $d_G(x,y)=n$ and $T_x$ a
breadth-first-search tree rooted at $x$ whose edges are all assigned
signature $\lambda=+1$. Then $T_x$ contains a $(x-y)$path $P$, a
$(x-y^{-})$path $P_1$, and a $(x-y^{+})$path $P_2$.

\vskip 0.2cm

\noindent{\bf Subcase 1.1}\quad Either $P_1\subset P$ or $P_2\subset
P$.

\vskip 0.2cm

 We may suppose that $P_1\subset P$ and $\lambda(y,y^{+})=-1$
without generality. If $P$ is not contained in $C$,then $P\cap
C=\{x=x_1,x_2,...,x_m=y\}$. By Lemma 2, we may choose an index $i$
such that
$C_i=x_i\overrightarrow{P}x_{i+1}\cup{x_i\overrightarrow{C}x_{i+1}}$
is an even $\Pi-$onesided cycle such that $P$ has no inner vertex on
$x_i\overrightarrow{C}x_{i+1}$ other than $x_i$ and $x_{i+1}$. Since
$|C_i|<|C|$, $C_i$ can't be a sum of two shorter cycles. Therefore,
$C_i$ satisfies (1) and further, $C_i$ is a fundamental cycle of a
breadth-first-tree rooted at some vertex of $C_i$. Now
$C_i\cup(P_2\cup P\cup \{(y,y^{+})\})$ has a shortest $\Pi-$twosided
cycle in $G$.

\vskip 0.2cm

If $P\subset C$, then $P_2$ can't be contained in $C$
 (since $C$ is $\Pi-$twosided). As we have reasoned above, there
 exists a segment $x_j\overleftarrow{C}x_{j+1}$ of $C$ such that
 $C_j=x_j\overleftarrow{C}x_{j+1}\cup{x_j\overleftarrow{P_2}x_{j+1}}$
 is a fundamental cycle and so, $C_j\cup(P_2\cup P\cup \{(y,y^{+})\})$ contains
 a shortest $\Pi-$twosided cycle in $G$.

\vskip 0.2cm

\noindent{\bf Subcase 1.2}\quad $P_1\not\subset P$ and
$P_2\not\subset P.$

\vskip 0.2cm

 If $\lambda(y,y^{-})=+1$ or
$\lambda(y,y^{+})=+1$, then $P_1\cup{P\cup {(y,y^{-})}} $ or
$P_2\cup{P\cup {(y,y^{+})}} $ will contain a shortest $\Pi-$twosided
cycle. Assume further that $\lambda(y,y^{-})=\lambda(y,y^{+})=-1 $.
Then $P_1\cup P_2\cup \{(y,y^{-}),(y,y^{+})\}$ also contains a
shortest $\Pi-$twosided cycle in $G$.

\vskip 0.2cm

\noindent{\bf Case 2.}\quad $|C|=2n+1,\,n\in N$

\vskip 0.2cm

 Let $x,y\in C$ with
$d_C(x,y)=d_C(x,y^{+})=n$ and $T_x$ be a breadth-first-search tree
rooted at $x$ with all its edges been assigned $\lambda=+1$. Then
$T_x$ contains a $(x-y)$path $P_1$ and a $(x,y^{+})$path $P_2$ and a
branched vertex $\alpha\in{P_1\cap P_2}$. Suppose that $P_1$ is not
a part of $C$ and $\lambda(y,y^{+})=-1$. Then as we have shown in
Case 1 there exists an index $i$ such that
$C_i=x_i\overrightarrow{C}x_{i+1}\cup{x_i\overrightarrow{P_1}x_{i+1}}$
is a fundamental cycle and thus $C_i\cup P_1\cup P_2\cup
\{(y,y^{+})\}$ contains a shortest $\Pi-$twosided cycle in $\cal C$.
This completes the proof of Lemma 3 and so, finishes the proof of
part (a) of Theorem A.

\vskip 0.2cm

As for the proof of parts (b) and (c) of Theorem A, it follows from
the proving procedure of (a) of Theorem A.

\vskip 0.2cm

Based on the above statement, Theorem A is proved.

\vskip 0.3cm

\vskip 0.5cm
\no{\bf References}
 \vskip 0.5cm
\def\hang{\hangindent\parindent}
\def\textindent#1{\indent\llap{#1\enspace}\ignorespaces}
\def\re{\par\hang\textindent}

\re{[1]}J.A.Bondy and U.S.R.Murty,\,{\it Graph theory with
applications},\,Macmillan,\,\\ London, (1978).
 \re{[2]}B.Mohar and
C.Thomassen,\,{\it Graphs on Surfaces}\, The Johns Hopkins
University Press,\,Baltimore and London,\, 2001.
 \re{[3]}C.Thomassen,\,Embeddings of graphs with no
short noncontractible cycles,\,{\it J. of Combin. Theory.,}\,Ser B
\bf\rm 48 (1990),\,pp155-177

\end{document}